\documentclass[11pt,reqno]{amsart}
\usepackage[T1]{fontenc}
\usepackage{lmodern}
\usepackage{amsmath,amssymb,amsthm,mathtools}
\usepackage{microtype}
\usepackage[hidelinks]{hyperref}
\hypersetup{pdftitle={Fermat quotients with CM Jacobians and trivial automorphism group},
 pdfauthor={Saeed Tafazolian},pdfsubject={Research manuscript},
 pdfkeywords={Fermat curves, automorphisms, complex multiplication, gonality, Jacobians, Schur polynomials}}
\usepackage[margin=1in]{geometry}
\allowdisplaybreaks[1]

\newtheorem{introthm}{Theorem}

\newtheorem{introcor}[introthm]{Corollary}

\newtheorem{theorem}{Theorem}[section]
\newtheorem{proposition}[theorem]{Proposition}
\newtheorem{lemma}[theorem]{Lemma}

\theoremstyle{remark}
\newtheorem{remark}[theorem]{Remark}

\newcommand{\Q}{\mathbb Q}
\newcommand{\C}{\mathbb C}
\newcommand{\F}{\mathbb F}
\newcommand{\Z}{\mathbb Z}
\newcommand{\Jac}{\operatorname{Jac}}
\newcommand{\Aut}{\operatorname{Aut}}
\newcommand{\End}{\operatorname{End}}
\newcommand{\Res}{\operatorname{Res}}
\newcommand{\Tr}{\operatorname{Tr}}
\newcommand{\sgn}{\operatorname{sgn}}
\newcommand{\col}{\operatorname{col}}
\newcommand{\pr}{\operatorname{pr}}
\newcommand{\gon}{\operatorname{gon}}
\newcommand{\lw}{\operatorname{lw}}
\newcommand{\doi}[1]{\href{https://doi.org/#1}{doi:#1}}

\title[Fermat quotients with trivial automorphism group]{Fermat quotients with CM Jacobians\\ and trivial automorphism group}
\author{Saeed Tafazolian}
\address{Institute of Mathematics, Statistics and Scientific Computing (IMECC),
University of Campinas (UNICAMP), Campinas, SP, Brazil}
\email{saeed@unicamp.br}
\date{}

\begin{document}

\begin{abstract}
Let $F_n$ be the Fermat curve of degree $n$ over $\mathbb Q$, and let
$D_n=F_n/S_3$ be its quotient by coordinate permutations. We prove that
$D_n$ has trivial geometric automorphism group for every $n\ge 11$
with $\gcd(n,6)=1$. Its Jacobian has complex multiplication.
For $n=11$, this Jacobian is absolutely simple of dimension five.
We construct an isogeny over $\mathbb Q$ of degree $11^5$ from
$\operatorname{Jac}(D_{11})$ to the Jacobian of $y^{11}=x(1-x)^2$.

For primes $p\equiv 1\pmod 3$, we also consider the metacyclic
quotients $C_p$ of $F_p$ over $\mathbb Q$. By a theorem of Carocca,
Lange and Rodr\'iguez, their Jacobians are absolutely simple and have
complex multiplication. We prove that $C_p$ has trivial geometric
automorphism group if and only if $p\ge 19$. Write $p=r^2+rs+s^2$
with $r,s$ positive integers. The natural degree-three cover
$Y_p\to C_p$ satisfies $\operatorname{gon}(Y_p)=r+s$, so
$\operatorname{gon}(C_p)\ge \lceil(r+s)/3\rceil$.
Since $g(C_p)=(p-1)/6$, both the genus and the gonality of $C_p$
tend to infinity with $p$.
\end{abstract}

\maketitle

\section{Introduction}

Moonen and Oort asked for a curve of genus at least four whose Jacobian
has complex multiplication (CM) and whose geometric automorphism group is trivial
\cite[Question~6.2]{MoonenOort}. See also \cite[Question~2.20]{Baldi}.
Bouchet, Hanselman, Pieper and Schiavone obtained infinitely many such
curves in genus four from CM points on a Shimura curve.
Their existence proof does not single out an individual CM fibre.
For the curves in their Table~2 they prove that the geometric
automorphism group is trivial. The evidence for complex multiplication
of these particular Jacobians is numerical
\cite[Remark~4.18 and Table~2]{BHPS}. In this paper we give explicit examples over
$\Q$ with absolutely simple Jacobians. Both the genus and the gonality
can be arbitrarily large.

Write $F_n$ for the smooth plane curve $X^n+Y^n+Z^n=0$. Its Jacobian
has CM. The same holds for the Jacobian of each of its quotients
\cite{KoblitzRohrlich}. The symmetric group $S_3$ acts by permuting
the coordinates. Let $D_n=F_n/S_3$ and let $\pi:F_n\to D_n$ be the
quotient map. Throughout, quotients denote smooth projective curves,
and gonality is taken over $\C$.

\begin{introthm}[Symmetric quotients]\label{thm:A}
For every integer $n\ge11$ with $\gcd(n,6)=1$, the curve $D_n$ is
defined over $\Q$, its Jacobian has CM, and
\[
 g(D_n)=\frac{(n-1)(n-5)}{12},\qquad \Aut_{\C}(D_n)=\{1\}.
\]
\end{introthm}

The first genera are $5,8,16,21$ and $33$. The degrees $5$ and $7$
give genus zero and one, respectively. For $n=11$ the rational
endomorphism algebra is $\End^0_{\C}(\Jac D_{11})=\Q(\zeta_{11})$.
It contains $\mu_{22}$, so its roots of unity alone do not exclude
automorphisms of order eleven. In Section~\ref{sec:eleven} we deduce
from Theorem~\ref{thm:A} and Torelli's theorem that $\zeta_{11}$
lies outside the integral endomorphism ring for the CM action specified there.

To prove Theorem~\ref{thm:A}, we use a theorem of Lim to show that
each automorphism of $D_n$ acts
diagonally in a basis of Schur differentials
\cite[Theorem~A]{Lim}. We need only four or five of these forms.
Up to a scalar, the chosen forms satisfy just one quadratic relation
of rank three. Its column space and the Schur coefficients force the
eigenvalues on the chosen forms to agree.
Their ratios generate $\C(D_n)$.

For $n>11$ prime to six, $\Jac D_n$ is not absolutely simple.
Indeed, the Fermat decomposition gives
factors of dimension $\varphi(n')/2$ for certain divisors $n'$ of $n$
\cite[pp.~1183--1184]{KoblitzRohrlich}. Every simple factor of
$\Jac F_n$ therefore has dimension at most $(n-1)/2$, whereas
$g(D_n)>(n-1)/2$.

For the second family let $p\equiv1\pmod3$ be prime. Let $k$ have
order three in $(\Z/p\Z)^\times$ and let $\zeta_p$ be a primitive
$p$th root of unity. On $F_p$ put
\[
 h(X:Y:Z)=(\zeta_p^{-k}X:\zeta_pY:Z),\qquad
 \tau(X:Y:Z)=(Y:Z:X).
\]
Put $G_p=\langle h,\tau\rangle$. This is a metacyclic group of order $3p$.
We write
\begin{equation}\label{eq:family}
 Y_p=F_p/\langle h\rangle,\qquad
 C_p=F_p/G_p=Y_p/\langle\bar\tau\rangle,\qquad
 K_p=\Q(\zeta_p)^{\langle k\rangle},
\end{equation}
where $\bar\tau$ is induced by $\tau$ and $k$ acts on $\zeta_p$ by
raising it to the $k$th power. The curves $C_p$ have genus
$(p-1)/6$. Carocca, Lange and Rodr\'iguez proved that $\Jac C_p$ is
absolutely simple with CM by $K_p$. They also proved that
$\Jac Y_p\sim_{\C}(\Jac C_p)^3$
\cite[Theorem~1, parts (1), (3) and~(4), p.~6160]{CLR}.

\begin{introthm}[Metacyclic quotients]\label{thm:B}
Let $p\equiv1\pmod3$ be prime. Write $p=r^2+rs+s^2$ with $r,s$ positive integers.
The curves $Y_p$ and $C_p$ in \eqref{eq:family} are defined over $\Q$, and
\begin{equation}\label{eq:main-gonality}
 \gon(Y_p)=r+s,\qquad
 \gon(C_p)\ge\left\lceil\frac{r+s}{3}\right\rceil.
\end{equation}
Moreover, $\Aut_{\C}(C_p)=\{1\}$ if and only if $p\ge19$.
\end{introthm}

The classical Hurwitz model recalled in \cite[p.~940]{TzermiasHK}
already gives $\gon(Y_p)\le r+s$. Indeed, put $m=r+s$ and $\ell=s$.
Then $m^2-m\ell+\ell^2=p$, and the equation $u^mv^\ell+v^m+u^\ell=0$ has
degree $m$ in $v$. Projection to $u$ gives the bound. We prove the
reverse inequality using the lattice-width theorem of Castryck and
Cools \cite{CastryckCools}. Since $(r+s)^2=p+rs$,
the bound gives $\gon(C_p)\ge3$ for $p>36$. The only roots of unity
in $K_p$ are $\pm1$. It follows that a nonhyperelliptic $C_p$ of
genus at least three has trivial automorphism group. We treat $p=19$
and $31$ by computing orders of canonical differentials.

\begin{introcor}\label{cor:C}
For all positive integers $g_0$ and $d_0$, there is a curve $C$ over
$\Q$ with $g(C)>g_0$, $\gon(C)>d_0$ and $\Aut_{\C}(C)=\{1\}$ whose
Jacobian is absolutely simple and has complex multiplication.
\end{introcor}

The genera of the metacyclic curves with trivial automorphism group
are exactly the integers $g\ge3$ for which $6g+1$ is prime. The
symmetric family gives other genera. For example, $D_{13}$ has genus
eight and trivial automorphism group, although its Jacobian is not simple.

\begin{remark}\label{rem:chebyshev}
Oort also asked for explicit CM curves with a small automorphism group
\cite[Question~11.2]{OortCM}. Let $H_d$ be the hyperelliptic curve
$y^2=(x+2)\varphi_d(x)$, where $\varphi_d(2\cos t)=2\cos(dt)$.
For $d=2^e\ge4$ or $d=\ell\ge5$ prime, its Jacobian is absolutely
simple with CM \cite[Theorem~1.1]{TTCM}. We show that the hyperelliptic
involution is the only nonidentity automorphism of $H_d$ over $\C$.

The finite automorphism group acts faithfully on the Jacobian, so it embeds
in the roots of unity of its endomorphism field. For $d=2^e$, this
field is $\Q(\zeta-\zeta^{-1})$, where $\zeta=\zeta_{4d}$.
The automorphism $\zeta\mapsto-\zeta^{-1}$ fixes this field and sends
$i$ to $-i$. Since all roots of unity in $\Q(\zeta)$ have order a
power of two, the only ones in the endomorphism field are $\pm1$.

For $d=\ell$ prime, the roots of unity in the endomorphism field
$\Q(\zeta_\ell)$ form $\mu_{2\ell}$. An extra automorphism would give a M\"obius transformation
of order $\ell$ on the hyperelliptic quotient. It permutes the
$\ell+1$ real branch points, namely $-2$ and the roots of $\varphi_\ell$.
It is thus defined over $\mathbb R$ and fixes a branch point since
$\ell\nmid\ell+1$. Sending that point to infinity gives a real affine
map, which cannot have odd finite order greater than one.
For $d>4$, the projective stabilizer of the roots of $\varphi_d$
alone is generated by $x\mapsto-x$ \cite[Remark~5.1]{TTSym}.
\end{remark}

The two families each contain a curve of genus five. Tzermias gave an
affine model for $D_{11}$. His proof also shows that $\Jac D_{11}$
is absolutely simple and has Mordell--Weil rank one over $\Q$
\cite[Lemma~4.1 and Remark~4.2]{Tzermias03}. The nine cyclic factors
of $\Jac F_{11}$ fall into two orbits under coordinate permutations
\cite[Proposition~2.1 and Lemma~2.5]{FGL}. In the notation of
\cite[\S2]{Tzermias03}, the three factors with $t=1,5,9$ have rank
zero. The six curves in the other orbit are $\Q$-isomorphic to the
smooth projective curve $Y_{11,2}$ given by $y^{11}=x(1-x)^2$.
This curve corresponds to $t=6$. Hence $\Jac D_{11}$ is $\Q$-isogenous
to $\Jac Y_{11,2}$. In Section~\ref{sec:eleven} we construct a specific
isogeny $\xi$ over $\Q$. For the canonical principal polarizations it
satisfies $\xi^\dagger\xi=[11]$, so $\deg\xi=11^5$. Over
$\Q(\zeta_{11})$ it respects the specified CM actions. The two
Jacobians have different integral endomorphism rings, so they are
not isomorphic over $\C$ even as unpolarized abelian varieties.
The curve $Y_{11,2}$ has an automorphism of order eleven, while
$D_{11}$ has trivial automorphism group.

\section{The symmetric Fermat quotient}\label{sec:symmetric}

In this section and the next, assume $n\ge5$ and $\gcd(n,6)=1$.
Put $F=F_n$ and let $\pi:F\to D_n=F/S_3$ be the quotient map. Every
automorphism of $F$ is a product of a diagonal automorphism and a
coordinate permutation \cite{Tzermias95}. The normalizer of $S_3$ in
this group is $S_3$ itself. Indeed, a diagonal element normalizing
$S_3$ commutes with each transposition, as one sees by projecting to
the permutation factor. Commuting with $(12)$ and $(23)$ forces the
three diagonal entries to agree. Thus every automorphism of $D_n$ that
lifts to $F$ is trivial. To study arbitrary automorphisms we use their
action on differentials.

\subsection{Invariant differentials and genus}
Let $\zeta=\zeta_n$ be a primitive $n$th root of unity. Write
$f=X^n+Y^n+Z^n$ and put
\[
 e_1=X+Y+Z,\qquad e_2=XY+YZ+ZX,\qquad e_3=XYZ,
\]
and
\[
 \Omega_0=X\,dY\wedge dZ-Y\,dX\wedge dZ+Z\,dX\wedge dY.
\]
For a homogeneous polynomial $P$ of degree $n-3$, put
$\omega_P=\Res_F(P\,\Omega_0/f)$. In the chart $Z\ne0$, with
$x=X/Z$ and $y=Y/Z$, this is the differential
$P(x,y,1)\,dx/(ny^{n-1})$, up to sign. By adjunction,
$P\mapsto\omega_P$ is an isomorphism from the homogeneous polynomials
of degree $n-3$ onto $H^0(F,\Omega_F^1)$. In particular,
\begin{equation}\label{eq:ratio}
 \frac{\omega_P}{\omega_{P'}}=\frac{P}{P'}
 \qquad\text{as rational functions on }F.
\end{equation}

For positive integers $a,b,c$ with $a+b+c=n$, write
$\omega_{a,b,c}=\omega_P$ with $P=X^{a-1}Y^{b-1}Z^{c-1}$. These forms
are a basis of $H^0(F,\Omega_F^1)$. Since the diagonal map
$(X:Y:Z)\mapsto(\zeta X:Y:Z)$ multiplies $\Omega_0$ by $\zeta$, it
acts on $\omega_{a,b,c}$ by $\zeta^a$. Similarly
$(X:Y:Z)\mapsto(X:\zeta Y:Z)$ acts by $\zeta^b$. A permutation $s$ of
the coordinates multiplies $\Omega_0$ by $\sgn(s)$, so
\begin{equation}\label{eq:perm}
 s^*\omega_{a,b,c}=\sgn(s)\,\omega_{s(a,b,c)}.
\end{equation}

By a \emph{triple} we mean $t=(a,b,c)$ with $a>b>c\ge1$. For a triple
with $a+b+c=n$, put
\[
 \eta_t=\sum_{s\in S_3}\sgn(s)\,\omega_{s(t)}.
\]
We identify $H^0(D_n,\Omega^1)$ with the $S_3$-invariant subspace of
$H^0(F,\Omega_F^1)$ via $\pi^*$.

\begin{lemma}\label{lem:basis}
The forms $\eta_t$, with $t$ running over the triples with
$a+b+c=n$, form a basis of $H^0(D_n,\Omega^1)$. Moreover
\[
 g(D_n)=\frac{(n-1)(n-5)}{12}.
\]
\end{lemma}

\begin{proof}
Let $\pr=\frac16\sum_{s\in S_3}s^*$ be the projection onto the
invariant subspace. By \eqref{eq:perm}, $\pr(\omega_{a,b,c})=0$ if two
of $a,b,c$ are equal, and $\pr(\omega_{s(t)})=\frac{\sgn(s)}6\,\eta_t$
for a triple $t$. Hence the $\eta_t$ span the invariant subspace. They
are linearly independent, since different triples involve disjoint
sets of basis vectors $\omega_{a,b,c}$.

The genus is a special case of Re's formulas for Fermat quotients
\cite{Re}. Here it follows by counting fixed points. A transposition,
say $(X:Y:Z)\mapsto(Y:X:Z)$, fixes the $n$ points of $F$ on the line
$X=Y$ and the point $(1:-1:0)$. A three-cycle fixes the two
points $(1:\zeta_3:\zeta_3^2)$ and $(1:\zeta_3^2:\zeta_3)$, where
$\zeta_3$ is a primitive cube root of unity. They lie on $F$ because
$3\nmid n$. Riemann--Hurwitz gives
\[
 n(n-3)=6\bigl(2g(D_n)-2\bigr)+3(n+1)+2\cdot2,
\]
and therefore $12g(D_n)=n^2-6n+5$.
\end{proof}

\subsection{Schur polynomials and the function field}
We describe $\eta_t$ by polynomials. For a triple $t=(a,b,c)$ put
\[
 A_t=\sum_{s\in S_3}\sgn(s)\,s\bigl(X^{a-1}Y^{b-1}Z^{c-1}\bigr)
 =\det\begin{pmatrix}
 X^{a-1}&X^{b-1}&X^{c-1}\\
 Y^{a-1}&Y^{b-1}&Y^{c-1}\\
 Z^{a-1}&Z^{b-1}&Z^{c-1}
 \end{pmatrix},
\]
and $\Delta=A_{(3,2,1)}=(X-Y)(X-Z)(Y-Z)$. The polynomial $A_t$ is
alternating, so it is divisible by $\Delta$, and
\[
 s_t=\frac{A_t}{\Delta}
\]
is a symmetric polynomial of degree $a+b+c-6$, the \emph{Schur
polynomial} of $t$. When $a+b+c=n$, the corresponding differential
satisfies
\begin{equation}\label{eq:eta-schur}
 \eta_t=\omega_{A_t}=\omega_{\Delta s_t}.
\end{equation}

\begin{lemma}\label{lem:schur}
Let $V_d$ be the space of symmetric polynomials of degree $d$.
\begin{enumerate}
\item The monomials $e_1^ie_2^je_3^l$ with $i+2j+3l=d$ form a basis
of $V_d$.
\item The Schur polynomials $s_t$ with $a+b+c=d+6$ form another basis
of $V_d$.
\item $s_{t+(1,1,1)}=e_3\,s_t$ for every triple $t$.
\end{enumerate}
\end{lemma}

\begin{proof}
Part (1) is the fundamental theorem on symmetric polynomials. For
(2), an alternating polynomial of degree $d+3$ is determined by its
coefficients at the monomials $X^{a-1}Y^{b-1}Z^{c-1}$ with $a>b>c$,
and $A_t$ is the only one of the polynomials $A_{t'}$ containing that
monomial. So the $A_t$ form a basis of the alternating polynomials of
degree $d+3$. Every such polynomial is $\Delta$ times a symmetric
polynomial, which gives (2). For (3), increasing all exponents by one
multiplies the rows of the determinant by $X,Y,Z$.
\end{proof}

\begin{proposition}\label{prop:model}
Put $u=e_2/e_1^2$ and $v=e_3/e_1^3$, and let $L=\Q$ or $\C$.
\begin{enumerate}
\item $L(D_n)=L(u,v)$.
\item Let $N_n\in\Z[u,v]$ be the polynomial such that
$x^n+y^n+z^n=N_n(u,v)$ whenever $x+y+z=1$, $xy+yz+zx=u$ and $xyz=v$.
Then $N_n$ is irreducible over $\C$, and $D_n$ is the smooth
projective model of $N_n(u,v)=0$.
\end{enumerate}
In particular, $D_n$ and the quotient map $\pi$ are defined over
$\Q$.
\end{proposition}

\begin{proof}
(1) Put $x=X/e_1$, $y=Y/e_1$ and $z=Z/e_1$. Then $L(F)=L(x,y,z)$,
the group $S_3$ permutes $x,y,z$, and $x,y,z$ are the roots of
$T^3-T^2+uT-v$. Hence $[L(F):L(u,v)]\le6$. Since $u,v$ are invariant
and $[L(F):L(F)^{S_3}]=6$, we get $L(D_n)=L(F)^{S_3}=L(u,v)$.

(2) The polynomial $N_n$ exists by Lemma~\ref{lem:schur}(1), applied
with $e_1=1$. Consider the map
\[
 \chi:\mathbb A^2\to\mathbb A^2,\qquad
 (x,y)\mapsto(u,v),
\]
computed with $z=1-x-y$. It is dominant, and
$N_n\circ\chi=x^n+y^n+(1-x-y)^n$, which is irreducible because $F$ is
smooth. If $N_n=GH$ with $G,H$ nonconstant, then
$N_n\circ\chi=(G\circ\chi)(H\circ\chi)$ with both factors
nonconstant, a contradiction. Since $u,v$ satisfy $N_n(u,v)=0$ on
$D_n$, part (1) shows that $D_n$ is the smooth projective model of
this curve.
\end{proof}

\subsection{Lim's theorem and Schur differentials}
For a finite morphism of curves $\nu:X\to X'$, the maps $\nu_*$ and
$\nu^*$ between Jacobians induce on holomorphic differentials,
respectively, the pullback $\nu^*$ and the trace $\Tr_\nu$. For an
automorphism $g$, the map $g_*$ acts on differentials as $g^*$.

\begin{theorem}[Lim]\label{thm:lim}
Let $m>3$ with $\gcd(m,6)=1$. Every element of
$\End^0(\Jac F_m)=\End(\Jac F_m)\otimes\Q$ is a $\Q$-linear
combination of the endomorphisms $g_*$, where $g$ runs over the
automorphisms of $F_m$.
\end{theorem}

This is the surjectivity of the map
$\Q[\Aut_{\C}(F_m)]\to\End^0_{\C}(\Jac F_m)$ in Lim's exact
sequence \cite[Theorem~A, p.~86]{Lim}. The surjectivity is proved
at the end of \cite[\S2, p.~94]{Lim}.

\begin{proposition}\label{prop:diagonal}
Let $\beta$ be an automorphism of $D_n$. For every triple $t$ with
$a+b+c=n$ there is $\lambda_t\in\C^\times$ such that
\[
 \beta^*\eta_t=\lambda_t\,\eta_t.
\]
\end{proposition}

\begin{proof}
Let $V_t\subset H^0(F,\Omega_F^1)$ be the span of the six forms
$\omega_{s(t)}$, $s\in S_3$. By \eqref{eq:perm} and the diagonal action on the monomial forms,
$V_t$ is stable under coordinate permutations and diagonal maps. They generate $\Aut(F)$ \cite{Tzermias95}, so
$g^*V_t\subseteq V_t$ for every $g\in\Aut(F)$. By the proof of
Lemma~\ref{lem:basis}, the only $S_3$-invariant forms in $V_t$ are the
multiples of $\eta_t$.

Consider the endomorphism
$\varepsilon=\pi^*\circ\beta_*\circ\pi_*$ of $\Jac F$. On differentials
it acts by $\omega\mapsto\pi^*\beta^*\Tr_\pi(\omega)$. Since
$\Tr_\pi(\pi^*\eta)=6\eta$ for $\eta\in H^0(D_n,\Omega^1)$, we get
\[
 \varepsilon^*(\pi^*\eta_t)=6\,\pi^*(\beta^*\eta_t),
\]
with $\eta_t$ viewed here as a differential on $D_n$. As $n>3$ and
$\gcd(n,6)=1$, Theorem~\ref{thm:lim} gives
$\varepsilon=\sum_g c_g\,g_*$ with $c_g\in\Q$. Thus
$\varepsilon^*=\sum_g c_g\,g^*$ maps $V_t$ into itself, and
$\beta^*\eta_t$ is an invariant form in $V_t$. Hence it is a multiple
of $\eta_t$, and the multiple is nonzero because $\beta^*$ is
invertible.
\end{proof}

\section{Automorphisms of the symmetric quotient}\label{sec:aut-sym}

\subsection{Quadratic relations in degrees four and five}
Let $d\in\{4,5\}$. Let $S_1,\dots,S_N$ be the monomial basis of $V_d$
from Lemma~\ref{lem:schur}(1), in the order
\[
 d=4:\ \ e_1^4,\ e_1^2e_2,\ e_2^2,\ e_1e_3;\qquad
 d=5:\ \ e_1^5,\ e_1^3e_2,\ e_1e_2^2,\ e_1^2e_3,\ e_2e_3.
\]
Let $U_d\subset V_d$ be the span of the monomials without $e_3$, that
is, of the first three. We identify $x\in\C^N$ with
$\sum x_iS_i\in V_d$.

A \emph{relation} is a symmetric $N\times N$ matrix $M=(m_{ij})$ such
that
\[
 \sum_{i,j}m_{ij}S_iS_j=0 .
\]
Its column space $\col(M)$ is then a subspace of $V_d$. In both
degrees, let $M_0$ be the matrix of the relation $S_1S_3-S_2^2$.

\begin{lemma}\label{lem:relations}
\begin{enumerate}
\item For $d=4$, every relation is a multiple of $M_0$.
\item For $d=5$, the relations are the matrices of
\[
 c_1(S_1S_3-S_2^2)+c_2(S_1S_5-S_2S_4)+c_3(S_2S_5-S_3S_4),
\]
and such a relation has rank at most three only if $c_2=c_3=0$.
\end{enumerate}
In both cases $M_0$ has rank three and $\col(M_0)=U_d$.
\end{lemma}

\begin{proof}
Each product $S_iS_j$ is a monomial in $e_1,e_2,e_3$, and distinct
monomials in these variables are linearly independent. So the
relations come from products that coincide. For $d=4$ the only
coincidence is $S_1S_3=S_2^2=e_1^4e_2^2$. For $d=5$ the three
coincidences are
\[
 S_1S_3=S_2^2=e_1^6e_2^2,\qquad
 S_1S_5=S_2S_4=e_1^5e_2e_3,\qquad
 S_2S_5=S_3S_4=e_1^3e_2^2e_3.
\]
Let $B$ be the matrix of the relation with coefficients
$c_1,c_2,c_3$. A direct expansion shows that its principal
$4\times4$ minor on the indices $2,3,4,5$ equals $c_3^4/16$, and the
one on $1,2,4,5$ equals $c_2^4/16$. So $B$ has rank at least four
unless $c_2=c_3=0$. Finally, in both cases $M_0$ is supported on the
first three rows and columns, where it is the nondegenerate matrix
\[
 \begin{pmatrix}0&0&\tfrac12\\0&-1&0\\\tfrac12&0&0\end{pmatrix}.
\]
Hence $\col(M_0)=U_d$.
\end{proof}

The triples $t$ with $a+b+c=d+6$ are
\[
 d=4:\ (7,2,1),(6,3,1),(5,4,1),(5,3,2);\qquad
 d=5:\ (8,2,1),(7,3,1),(6,4,1),(6,3,2),(5,4,2).
\]

\begin{lemma}\label{lem:table}
The coefficient matrices of the Schur polynomials in the monomial
bases $(S_1,\dots,S_N)$ are
\[
 \begin{pmatrix}
 1&-3&1&2\\
 0&1&-1&-1\\
 0&0&1&-1\\
 0&0&0&1
 \end{pmatrix}
 \qquad\text{and}\qquad
 \begin{pmatrix}
 1&-4&3&3&-2\\
 0&1&-2&-1&2\\
 0&0&1&-1&-1\\
 0&0&0&1&-1\\
 0&0&0&0&1
 \end{pmatrix}
\]
for degrees four and five respectively. The rows correspond to the
triples in the order just given.
\end{lemma}

\begin{proof}
Expand $A_t=\Delta s_t$, or apply the Jacobi--Trudi formula
\cite[Chapter~I, (3.4)]{Macdonald}. For example,
$s_{(5,3,2)}=e_3\,s_{(4,2,1)}=e_3e_1$ by Lemma~\ref{lem:schur}(3).
\end{proof}

\subsection{The action on the function field}

\begin{proof}[Proof of Theorem~\ref{thm:A}]
The genus is given by Lemma~\ref{lem:basis}, and the field of
definition by Proposition~\ref{prop:model}. The Jacobian of $F_n$ is
isogenous to a product of abelian varieties with CM
\cite{KoblitzRohrlich}. Since $\pi_*\circ\pi^*$ is multiplication by
six on $\Jac D_n$, the map $\pi^*$ has finite kernel, and $\Jac D_n$
is isogenous to an abelian subvariety of $\Jac F_n$. Hence $\Jac D_n$
also has CM.

Now let $n\ge11$ and let $\beta$ be an automorphism of $D_n$. Put
\[
 (d,r)=\begin{cases}
 (4,\,(n-10)/3) & \text{if } n\equiv1\pmod3,\\
 (5,\,(n-11)/3) & \text{if } n\equiv2\pmod3.
 \end{cases}
\]
Then $r\ge0$ is an integer, $d+3r=n-6$, and $2d<n$. At $n=11$ the last inequality reads $10<11$. For
$S\in V_d$, define the differential $\langle S\rangle$ on $D_n$ by
\[
 \pi^*\langle S\rangle=\omega_{\Delta e_3^rS}.
\]
By \eqref{eq:eta-schur} and Lemma~\ref{lem:schur}(3), for each triple
$t$ with $a+b+c=d+6$ we have
\[
 \langle s_t\rangle=\eta_{t+(r,r,r)} .
\]
By Proposition~\ref{prop:diagonal}, $\beta^*$ maps the space
$W=\{\langle S\rangle:S\in V_d\}$ into itself. Since
$S\mapsto\langle S\rangle$ is injective by adjunction, the action on
$W$ defines an invertible map $\psi:V_d\to V_d$ by
$\beta^*\langle S\rangle=\langle\psi(S)\rangle$. Then
\begin{equation}\label{eq:psi-diagonal}
 \psi(s_t)=\lambda_t\,s_t\qquad(a+b+c=d+6),
\end{equation}
where we write $\lambda_t$ for $\lambda_{t+(r,r,r)}$.

\smallskip
\noindent\emph{Step 1. $\psi(U_d)=U_d$.}
Let $\Psi=(\psi_{ij})$ be the matrix of $\psi$ in the basis
$S_1,\dots,S_N$, so that $\psi(S_j)=\sum_i\psi_{ij}S_i$. We claim that
$\Psi M\Psi^T$ is a relation whenever $M$ is. Let $M$ be a relation.
By \eqref{eq:ratio},
\[
 \frac{\langle S_i\rangle\langle S_j\rangle}{\langle e_1^d\rangle^2}
 =\frac{S_iS_j}{e_1^{2d}},
\]
so the quadratic differential
$\sum m_{ij}\langle S_i\rangle\langle S_j\rangle$ vanishes. Applying
$\beta^*$, which is multiplicative, we obtain
\[
 \sum_{i,j}m_{ij}\langle\psi(S_i)\rangle\langle\psi(S_j)\rangle=0 .
\]
Dividing by $\langle e_1^d\rangle^2$ again, the homogeneous
polynomial $\sum m_{ij}\psi(S_i)\psi(S_j)$ of degree $2d$ vanishes on
$F$. Since $f$ divides it and $\deg f=n>2d$, the polynomial is
zero. Expanding $\psi(S_i)$ shows that this polynomial is
$\sum_{i',j'}(\Psi M\Psi^T)_{i'j'}S_{i'}S_{j'}$, as claimed.

Apply this to the relation $M_0$ of Lemma~\ref{lem:relations}. The
relation $\Psi M_0\Psi^T$ has rank three, so it equals $\kappa M_0$
for some $\kappa\ne0$, by Lemma~\ref{lem:relations}. Comparing column
spaces, $\Psi(U_d)=U_d$.

\smallskip
\noindent\emph{Step 2. $\psi$ is a scalar.}
By Step~1, $\psi$ induces an automorphism of $V_d/U_d$. The images of
the Schur polynomials in this quotient are eigenvectors, with the same
eigenvalues as in \eqref{eq:psi-diagonal}. For $d=4$ the quotient is
one-dimensional, and the last column of the first matrix in
Lemma~\ref{lem:table} shows that all four images are nonzero. For
$d=5$ the quotient is two-dimensional. The last two columns of the
second matrix give the five nonzero vectors
\[
 (3,-2),\quad(-1,2),\quad(-1,-1),\quad(1,-1),\quad(0,1).
\]
They span five distinct lines. A linear map on a two-dimensional space
preserving three distinct lines is scalar. In both cases the induced
automorphism is therefore scalar, and every $\lambda_t$ is its scalar
eigenvalue. Since the $s_t$ form a basis of $V_d$, the map $\psi$
itself is scalar.

\smallskip
\noindent\emph{Step 3. $\beta=1$.}
Since $d\ge4$, the space $V_d$ contains $e_1^d$, $e_1^{d-2}e_2$ and
$e_1^{d-3}e_3$. By \eqref{eq:ratio},
\[
 u=\frac{e_2}{e_1^2}=\frac{\langle e_1^{d-2}e_2\rangle}{\langle e_1^d\rangle},
 \qquad
 v=\frac{e_3}{e_1^3}=\frac{\langle e_1^{d-3}e_3\rangle}{\langle e_1^d\rangle}.
\]
By Step~2, $\beta^*u=u$ and $\beta^*v=v$. By
Proposition~\ref{prop:model}(1), $u$ and $v$ generate $\C(D_n)$.
Hence $\beta$ is the identity.
\end{proof}

\section{The genus-five curve \texorpdfstring{$D_{11}$}{D11}}\label{sec:eleven}

Here $F=F_{11}$, $\pi:F\to D_{11}$, $\zeta=\zeta_{11}$ and
$K=\Q(\zeta)$. By Lemma~\ref{lem:basis}, $g(D_{11})=5$, and a basis of
$H^0(D_{11},\Omega^1)$ is given by $\eta_t$ for the five triples
\begin{equation}\label{eq:triples}
 (8,2,1),\quad(7,3,1),\quad(6,4,1),\quad(6,3,2),\quad(5,4,2).
\end{equation}

For the Fermat triple $(1,2,8)$ the classical CM type is
$H_{1,2,8}=\{m:\langle m\rangle+\langle2m\rangle+\langle8m\rangle=11\}$,
where $\langle a\rangle$ is the least positive residue modulo $11$ and
$m\in(\Z/11\Z)^\times$ \cite[\S1, pp.~1183--1184]{KoblitzRohrlich}.
Thus $H_{1,2,8}=\{1,2,3,6,7\}$. We identify the action on
$\Jac D_{11}$ that gives this type.

\begin{proposition}\label{prop:CM}
Put $\theta=\frac13(\zeta+\zeta^2+\zeta^8)$. Then
$\End^0_{\C}(\Jac D_{11})\cong K$ and $\Jac D_{11}$ is absolutely
simple. With respect to the embedding of $K$ that sends $\theta$ to
the endomorphism $\Theta$ below, the CM type of $\Jac D_{11}$ is
\[
 \Phi=\{1,2,3,6,7\}\subset(\Z/11\Z)^\times,
\]
where $m$ denotes the embedding $\sigma_m:\zeta\mapsto\zeta^m$.
\end{proposition}

\begin{proof}
Let $\rho(X:Y:Z)=(\zeta X:Y:Z)$ and
$\Theta=\frac16\,\pi_*\circ\rho_*\circ\pi^*\in\End^0_{\C}(\Jac D_{11})$.
On differentials, $\pi^*(\Theta^*\eta)=\pr(\rho^*(\pi^*\eta))$, where
$\pr$ is the projection in the proof of Lemma~\ref{lem:basis}.
Since each of $a,b,c$ occurs twice as the first entry of $s(a,b,c)$,
we get
\[
 \Theta^*\eta_{(a,b,c)}=\frac{\zeta^a+\zeta^b+\zeta^c}{3}\,\eta_{(a,b,c)} .
\]
Multiplying $\{1,2,8\}$ by $m\in\Phi$ modulo $11$ gives exactly the
five sets in \eqref{eq:triples}. So the eigenvalues of $\Theta^*$ on
$H^0(D_{11},\Omega^1)$ are $\sigma_m(\theta)$ for $m\in\Phi$.

The element $\theta$ generates $K$. Indeed, if
$\sigma_m(\theta)=\sigma_j(\theta)$, then, since the only rational
relations among $1,\zeta,\dots,\zeta^{10}$ are multiples of
$1+\zeta+\cdots+\zeta^{10}$, we get $\{m,2m,8m\}=\{j,2j,8j\}$ in
$\F_{11}$. The set $\{1,2,8\}$ has trivial stabilizer under
multiplication by $\F_{11}^\times$. Indeed, a stabilizing element must
belong to $\{1,2,8\}$, while $2\{1,2,8\}=\{2,4,5\}$ and
$8\{1,2,8\}=\{5,8,9\}$. Hence $m=j$, and $[\Q(\theta):\Q]=10$.

The eigenvalues of $\Theta$ on $H^1(D_{11},\C)=H^{1,0}\oplus H^{0,1}$
are the $\sigma_m(\theta)$ for $m\in\Phi$ and their complex conjugates
$\sigma_{-m}(\theta)$. Since $\Phi\sqcup(-\Phi)=(\Z/11\Z)^\times$,
these are the ten conjugates of $\theta$. So the characteristic
polynomial of $\Theta$ on $H^1(D_{11},\Q)$ is the minimal polynomial of
$\theta$, and $\Q[\Theta]\cong K$. This gives the embedding of $K$ and
the CM type $\Phi$.

Since $K/\Q$ is abelian, primitivity of the CM type is equivalent
to triviality of its stabilizer. If $j\Phi=\Phi$ with
$j\in(\Z/11\Z)^\times$, summing the elements gives $8j\equiv8$, so
$j=1$. Hence $\Jac D_{11}$ is absolutely simple
\cite[\S1]{KoblitzRohrlich}. Finally, $\End^0_{\C}(\Jac D_{11})$ is a
division algebra containing $K$ and acting faithfully on the
ten-dimensional space $H^1(D_{11},\Q)$, so its dimension over $\Q$
divides ten. As $[K:\Q]=10$, it equals $K$.
\end{proof}

\begin{proposition}\label{prop:Y}
Let $Y_{11,2}$ be the smooth projective model of $y^{11}=x(1-x)^2$.
Then $Y_{11,2}$ has genus five and $\Jac D_{11}$ is isogenous over
$\Q$ to $\Jac Y_{11,2}$. The curves $D_{11}$ and $Y_{11,2}$ are not
isomorphic over $\C$.
\end{proposition}

\begin{proof}
The function $x$ gives a cyclic cover of degree eleven over $\C$,
totally ramified at $0,1,\infty$ and unramified elsewhere, since
$x(1-x)^2$ has orders $1,2,-3$ at these points. Thus
$2g(Y_{11,2})-2=-22+3\cdot10=8$.

On the chart $Z\ne0$ of $F$, put $a=X/Z$ and $b=Y/Z$. The formulas
\[
 q:F\longrightarrow Y_{11,2},\qquad (x,y)=(-a^{11},-ab^2)
\]
satisfy the equation of $Y_{11,2}$, since $1-x=-b^{11}$. The resulting
rational map extends to a morphism of smooth projective curves over
$\Q$. Over $\C$ it is the quotient by the group generated by
$h_2(X:Y:Z)=(\zeta^{-2}X:\zeta Y:Z)$. Indeed, $h_2$ fixes $x,y$ and has
order eleven. Adjoining $b$ recovers $a=-y/b^2$, with $b^{11}=x-1$, so
the function field extension has degree eleven.

Consider $\xi=q_*\pi^*:\Jac D_{11}\to\Jac Y_{11,2}$, defined over $\Q$.
On $\omega_{1,2,8}$ the eigenvalue of $h_2^*$ is $\zeta^{-2}\zeta^2=1$.
Consequently $\omega_{1,2,8}=q^*\omega_Y$ for a holomorphic
differential $\omega_Y$ on $Y_{11,2}$. Since
$\xi^*\omega_Y=\Tr_\pi(q^*\omega_Y)$, we have
\[
 \pi^*(\xi^*\omega_Y)=\sum_{s\in S_3}s^*\omega_{1,2,8}
 =-\pi^*\eta_{(8,2,1)}\ne0.
\]
Thus $\xi$ is nonzero. Since $\Jac D_{11}$ is absolutely simple by
Proposition~\ref{prop:CM} and both Jacobians have dimension five,
$\xi$ is an isogeny.

The automorphism $\alpha(x,y)=(x,\zeta y)$ of $Y_{11,2}$ has order
eleven. Since $\Aut_{\C}(D_{11})=\{1\}$ by Theorem~\ref{thm:A}, the two
curves are not isomorphic over $\C$.
\end{proof}

The cyclic equation also gives the CM type directly. Write
$P_0,P_1,P_\infty$ for the respective points of $Y_{11,2}$ above
$0,1,\infty$. The orders of $y$ at $P_0,P_1,P_\infty$ are $1,2,-3$.
For $\omega=(1-x)^i\,dx/y^l$ one finds
\[
 \bigl(\operatorname{ord}_{P_0},\operatorname{ord}_{P_1},
 \operatorname{ord}_{P_\infty}\bigr)(\omega)
 =(10-l,\ 11i+10-2l,\ 3l-11i-12).
\]
Away from these three points such a form is regular, since $x$ is
unramified there and $y$ has no zeros or poles. Hence the five forms
\[
 \frac{dx}{y^4},\quad\frac{dx}{y^5},\quad\frac{(1-x)\,dx}{y^8},\quad
 \frac{(1-x)\,dx}{y^9},\quad\frac{(1-x)\,dx}{y^{10}}
\]
are holomorphic, with orders $(6,2,0)$, $(5,0,3)$, $(2,5,1)$,
$(1,3,4)$ and $(0,1,7)$ at $(P_0,P_1,P_\infty)$. They are eigenvectors
of $\alpha^*$ with eigenvalues $\zeta^{-l}=\zeta^m$ for
$m=7,6,3,2,1$, so they form a basis of $H^0(Y_{11,2},\Omega^1)$. The
quotient of $Y_{11,2}$ by $\langle\alpha\rangle$ is rational, so
$1+\alpha_*+\cdots+\alpha_*^{10}=0$ on $\Jac Y_{11,2}$. Since
$\alpha_*\ne1$, the map $\zeta\mapsto\alpha_*$ embeds $K$ in
$\End^0_{\C}(\Jac Y_{11,2})$, with CM type $\Phi$. By
Proposition~\ref{prop:Y}, $\End^0_{\C}(\Jac Y_{11,2})$ is a field of
degree ten, so it equals $\Q(\alpha_*)$.

\begin{proposition}\label{prop:equivariant}
The isogeny $\xi$ of Proposition~\ref{prop:Y} is defined over $\Q$.
After base change to $K$, it respects the two CM actions
$\theta\mapsto\Theta$ on $\Jac D_{11}$ and
$\zeta\mapsto\alpha_*$ on $\Jac Y_{11,2}$. In particular,
$\xi\circ\Theta=\theta(\alpha_*)\circ\xi$ in
$\operatorname{Hom}^0(\Jac D_{11},\Jac Y_{11,2})$, where
$\theta(\alpha_*)=\frac13(\alpha_*+\alpha_*^2+\alpha_*^8)$.
With respect to the canonical principal polarizations, the adjoint
$\xi^\dagger$ satisfies
\begin{equation}\label{eq:xi-degree}
 \xi^\dagger\xi=[11],\qquad \deg\xi=11^5.
\end{equation}
Its geometric kernel is a maximal isotropic subgroup of
$(\Jac D_{11})[11]$ for the Weil pairing. It is isomorphic to
$(\Z/11\Z)^5$.
\end{proposition}

\begin{proof}
Let $\iota_D:K\to\End^0_{\C}(\Jac D_{11})$ and
$\iota_Y:K\to\End^0_{\C}(\Jac Y_{11,2})$ be the embeddings with
$\iota_D(\theta)=\Theta$ and $\iota_Y(\zeta)=\alpha_*$. Both embeddings are isomorphisms. Conjugation by $\xi$ gives an
automorphism $\sigma_c$ of $K$ with $\xi\circ\iota_D(x)=\iota_Y(\sigma_c(x))\circ\xi$ for all
$x\in K$. The isomorphism $\xi^*:H^0(Y_{11,2},\Omega^1)\to
H^0(D_{11},\Omega^1)$ therefore identifies the eigenvalues of
$\iota_Y(\zeta^c)^*$ with those of $\iota_D(\zeta)^*$. Both CM types are
$\Phi$, so these eigenvalues are $\zeta^{cm}$ and $\zeta^m$ for
$m\in\Phi$, respectively. Hence $c\Phi=\Phi$, and $c=1$ because the
stabilizer of $\Phi$ is trivial (proof of Proposition~\ref{prop:CM}).

Put $H=\langle h_2\rangle$. The norm is adjoint to pullback, so
$\xi^\dagger=\pi_*q^*$. For each basis form $\eta_t$ the trace formula gives
\[
 \pi^*(\xi^\dagger\xi)^*\eta_t
 =\left(\sum_{s\in S_3}s^*\right)
  \left(\sum_{h\in H}h^*\right)\pi^*\eta_t .
\]
The $H$-trace kills $\omega_{a,b,c}$ unless $b\equiv2a\pmod{11}$,
in which case it multiplies the form by eleven. The six summands of
$\pi^*\eta_t$ are indexed by the orderings of $m\{1,2,8\}$ for some
$m\in\Phi$. Of $1,2,8$, only $1$ has its double in the set. Thus the
$H$-trace of $\pi^*\eta_t$ is $11\nu$, where $\nu$ is one of its
signed summands. By \eqref{eq:perm}, $\sum_{s\in S_3}s^*\nu=\pi^*\eta_t$.
It follows that $(\xi^\dagger\xi)^*\eta_t=11\eta_t$ for all $t$.
In characteristic zero this gives $\xi^\dagger\xi=[11]$. Since
$\deg\xi^\dagger=\deg\xi$ and $g(D_{11})=5$, we get
$(\deg\xi)^2=11^{10}$.

Write $A=\Jac D_{11}$ and $B=\Jac Y_{11,2}$. The group
$H_\xi=\ker(\xi)(\C)$ lies in $A[11]$ and has order $11^5$.
Since $(\xi\xi^\dagger)\xi=\xi[11]_A=[11]_B\xi$ and $\xi$ is
surjective, $\xi\xi^\dagger=[11]_B$. Hence
$\xi^\dagger(B[11])\subseteq H_\xi$ and $\ker(\xi^\dagger)\subseteq B[11]$.
As $\deg\xi^\dagger=11^5$, the two groups in the first inclusion
have the same order. For $P\in B[11]$ and $Q\in H_\xi$,
functoriality of the Weil pairing gives
$e_{11,A}(\xi^\dagger P,Q)=e_{11,B}(P,\xi Q)=1$.
Thus $H_\xi$ is isotropic of dimension five in the ten-dimensional
symplectic space $A[11]$ over $\F_{11}$.
\end{proof}

The two Jacobians also have different integral endomorphism rings
over $\C$. On $\Jac Y_{11,2}$, the automorphism $\alpha$ gives the
full ring $\Z[\zeta]=\mathcal O_K$. On $\Jac D_{11}$, put
$z_D=\iota_D(\zeta)$, with $\iota_D$ as in the preceding proof.
Since $\xi z_D=\alpha_*\xi$ and $\alpha_*$ preserves the canonical
polarization, \eqref{eq:xi-degree} gives $z_D^\dagger z_D=1$.
If $z_D$ belonged to $\End_{\C}(\Jac D_{11})$, it would be a polarized
automorphism of order eleven. By Torelli's theorem and $\Aut_{\C}(D_{11})=\{1\}$,
the only polarized automorphisms are $\pm1$. Hence
$\End_{\C}(\Jac D_{11})$ is a proper order in $K$.

\section{Metacyclic Fermat quotients}\label{sec:metacyclic}

We return to the quotients $Y_p$ and $C_p$ of \eqref{eq:family}. In
Sections~\ref{sec:metacyclic}--\ref{sec:plane}, $p\equiv1\pmod3$ is
prime, $\zeta=\zeta_p$, and $1\le k\le p-2$ represents one of the two
elements of order three in $(\Z/p\Z)^\times$. In particular,
$k^2+k+1\equiv0\pmod p$.

\begin{proposition}\label{prop:metacyclic}
The group $G_p=\langle h,\tau\rangle$ is a semidirect product of order
$3p$, with $\tau h\tau^{-1}=h^k$. A model for $Y_p$ is
\begin{equation}\label{eq:kummer}
 y^p=x(1-x)^k .
\end{equation}
The map $x$ has degree $p$ and is totally ramified over $0,1,\infty$.
The automorphism $\bar\tau$ of $Y_p$ has two fixed points, and
\[
 g(Y_p)=\frac{p-1}{2},\qquad g(C_p)=\frac{p-1}{6}.
\]
\end{proposition}

\begin{proof}
Conjugation by $\tau$ permutes the diagonal entries of $h$ to
$(\zeta,1,\zeta^{-k})$. Scaling the entries by $\zeta^k$ gives $h^k$,
since $k+1\equiv-k^2\pmod p$. On the chart $Z\ne0$, put $a=X/Z$ and
$b=Y/Z$. The functions
\[
 x=-a^p,\qquad y=(-1)^{k+1}ab^k
\]
are $h$-invariant and satisfy \eqref{eq:kummer}. Adjoining $b$
recovers $a$, and $b^p=x-1$. Consequently $[\C(F_p):\C(x,y)]\le p$,
with equality since $h$ has order $p$. This proves the model. The
orders of $x(1-x)^k$ at $0,1,\infty$ are $1,k,-k-1$, all prime to $p$.
These are the only branch values of $x$, with ramification index $p$
at each, so Riemann--Hurwitz gives $2g(Y_p)-2=-2p+3(p-1)=p-3$.

On the $x$-line, $\bar\tau$ acts by $x\mapsto(x-1)/x$. Its two fixed
points are outside $\{0,1,\infty\}$. Write $\alpha(x,y)=(x,\zeta y)$
for the deck generator. It is induced by
$\rho(X:Y:Z)=(\zeta X:Y:Z)$, and $\tau\rho\tau^{-1}=\rho^{k^2}h^{-1}$
in projective coordinates. Thus
$\bar\tau\alpha\bar\tau^{-1}=\alpha^{k^2}$. On each of the two fixed
fibres, after choosing a point to identify the fibre with $\F_p$, the
action is $j\mapsto k^2j+b_0$ for some $b_0\in\F_p$. It has one fixed
point since $k^2\ne1$. Riemann--Hurwitz for $Y_p\to C_p$ now gives
$p-3=3(2g(C_p)-2)+4$.
\end{proof}

In the notation of Proposition~\ref{prop:Y}, $Y_p$ is the curve
$Y_{p,k}:y^p=x(1-x)^k$. Both families thus involve cyclic quotients of
Fermat curves of this shape.

\begin{remark}\label{rem:descent}
Fit\'e, Gonz\'alez and Lario describe the order-three quotient over
$\Q$ \cite[Remark~2.8]{FGL}. In our model both quotient maps are also
defined over $\Q$. Indeed, over
$\Q(\zeta_p)$ every Galois automorphism sends $h$ to a power of $h$
and fixes $\tau$, so the two invariant function fields are stable
under Galois descent. Interchanging $X$ and $Y$ sends the exponent
pair $(-k,1)$ to $(1,-k)=-k(-k^{-1},1)$ and conjugates $\tau$ to its
inverse. The two choices of $k$ therefore give isomorphic quotients
over $\Q$.
\end{remark}

\begin{proposition}[Carocca--Lange--Rodr\'iguez]\label{prop:knownCM}
The Jacobian $\mathcal J_p=\Jac(C_p)$ is absolutely simple,
$\End^0_{\C}(\mathcal J_p)=K_p$, and
\begin{equation}\label{eq:cube}
 \Jac(Y_p)\sim_{\C}\mathcal J_p^3.
\end{equation}
\end{proposition}

\begin{proof}
On $Y_p$ the automorphisms $\alpha$ and $\bar\tau$ generate a
metacyclic group of order $3p$. The quotient is the $x$-line modulo
$x\mapsto(x-1)/x$. The orbit $\{0,1,\infty\}$ and the two fixed values
of this transformation give the signature $(p,3,3)$. Apply
\cite[Theorem~1, p.~6160]{CLR} with $q=p$, $n=3$, $Y=Y_p$ and $X=C_p$.
Part~(1) gives CM by $K_p$, and part~(4) gives absolute simplicity.
Part~(3) gives \eqref{eq:cube}.

Put $D=\End^0_{\C}(\mathcal J_p)$. Since $\mathcal J_p$ is simple,
$D$ is a division algebra. The space $H_1(\mathcal J_p,\Q)$ is a
$D$-module, so $\dim_{\Q}D$ divides $2\dim\mathcal J_p=(p-1)/3$.
But $K_p\subseteq D$ and $[K_p:\Q]=(p-1)/3$. Hence $D=K_p$.
\end{proof}

Hidalgo and Rodr\'iguez give the order-three automorphism on
$y^p=x^\gamma(x-1)$ for $\gamma^2+\gamma+1\equiv0\pmod p$.
They record its two fixed points and the genus $(p-1)/6$ of the
quotient \cite[Remark~2(2), p.~335]{HidalgoRodriguez}. The change
$(x,y)\mapsto(1-x,-y)$ gives our model \eqref{eq:kummer} with
$k=\gamma$. They also state \eqref{eq:cube} in
\cite[\S4.1, p.~339]{HidalgoRodriguez}.

In the notation of Tzermias, $Y_p$ is the Hurwitz--Klein curve
$F_{p,k}$ after the change $x\mapsto1-x$, and $C_p$ is its quotient
$C_0$ \cite[pp.~939--941]{TzermiasHK}. Over $\Q(\zeta_p)$, Lim proved
that $\Jac Y_p$ is isomorphic to $(\Jac C_p)^3$ as an abelian variety
\cite[Corollary~6.3]{LimCyclic}. Tzermias describes this isomorphism
in \cite[\S2, p.~941]{TzermiasHK}.

\section{Gonality and automorphisms of the metacyclic quotients}\label{sec:gonality}

\begin{lemma}\label{lem:roots-proper}
The only roots of unity in $K_p$ are $1$ and $-1$.
\end{lemma}

\begin{proof}
The roots of unity in $\Q(\zeta_p)$ form $\mu_{2p}$. Every member
other than $\pm1$ generates $\Q(\zeta_p)$ and cannot belong to its
proper subfield $K_p$.
\end{proof}

\begin{lemma}\label{lem:aut-field}
Let $X$ be a nonhyperelliptic complex curve of genus at least three. If
$\End^0_{\C}(\Jac X)$ is a field whose only roots of unity are $\pm1$,
then $\Aut_{\C}(X)=\{1\}$.
\end{lemma}

\begin{proof}
Each automorphism has finite order and acts as $1$ or $-1$ on the
Jacobian. It therefore acts trivially on the projective canonical
space. The canonical map is an equivariant embedding.
\end{proof}

\begin{proposition}\label{prop:toric-gonality}
The gonality formula and bound in \eqref{eq:main-gonality} hold.
In particular, $\Aut_{\C}(C_p)=\{1\}$ for $p>36$.
\end{proposition}

\begin{proof}
Write $p=r^2+rs+s^2$ in positive integers $r,s$, and put $m=r+s$, $\ell=s$.
Such a representation exists since $p$ splits in the Eisenstein
integers \cite[p.~940]{TzermiasHK}. Then $p=m^2-m\ell+\ell^2$ and
$\gcd(m,\ell)=1$. By Remark~\ref{rem:descent} we may choose
$k\equiv-m/\ell\pmod p$, which satisfies $k^2+k+1\equiv0\pmod p$.
Put $t=1-x$ in \eqref{eq:kummer}. The Newton triangle of
$1-t-t^{-k}y^p$ has vertices $(0,0),(1,0),(-k,p)$.
It is equivalent to the triangle $\Pi$ with vertices
$(\ell,0),(m,\ell),(0,m)$ of the Hurwitz model. Indeed, put
$u=(m-\ell,\ell)$ and $v=(-\ell,m)$. Then $v+ku\in p\Z^2$ and
$\det(u,v)=p$. The matrix with columns $u,(v+ku)/p$ is unimodular
and sends the original triangle to $\Pi-(\ell,0)$.

The trinomial is nondegenerate since its $y$-derivative never vanishes
on the torus and the edge polynomials are binomials with primitive
exponent differences. By Pick's theorem, $\Pi$ has $(p-1)/2$ interior
lattice points. Thus their convex hull cannot be equivalent to
$\Upsilon$, which has four lattice points. The primitive edges also
exclude $d\Sigma$ for $d\ge2$. By \cite[Lemma~5.2(i) and
Corollary~6.2]{CastryckCools}, $\gon(Y_p)=\lw(\Pi)$.

Let $(a,b)$ be a primitive direction. After subtracting $\ell a$, the
values of its linear form at the three vertices are $0,D,E$, where
$D=(m-\ell)a+\ell b$ and $E=-\ell a+mb$. Let $w$ be their range. Then
\[
 p(a^2-ab+b^2)=D^2-DE+E^2\le w^2.
\]
Indeed, if $h,j$ are the successive gaps between $0,D,E$ in increasing
order, the middle expression is $h^2+hj+j^2\le(h+j)^2$.
The directions with $a^2-ab+b^2=1$ are, up to sign, $(1,0),(0,1)$
and $(1,1)$. Each gives width $m$. Since $a^2-ab+b^2\equiv(a+b)^2\pmod3$, the norm cannot be two.
Every other direction therefore has $w^2\ge3p>m^2$, using
$4p-3m^2=(m-2\ell)^2$. Hence $\gon(Y_p)=m=r+s$.

The map $Y_p\to C_p$ has degree three. Therefore $3\gon(C_p)\ge r+s$.
For $p>36$ we have $r+s>\sqrt p>6$, so $C_p$ is nonhyperelliptic.
Proposition~\ref{prop:knownCM} and Lemmas~\ref{lem:roots-proper}
and~\ref{lem:aut-field} give $\Aut_{\C}(C_p)=\{1\}$.
\end{proof}

For $p=7,13$ the formula gives the known values $\gon(Y_7)=3$ and
$\gon(Y_{13})=4$ \cite[p.~946]{TzermiasHK}. For every admissible
$p\ge13$ it gives $\gon(Y_p)\ge4$, the gonality hypothesis in
\cite[Lemma~3.1(b)]{TzermiasHK}.

\section{The small primes}\label{sec:plane}

\begin{lemma}\label{lem:small}
The curves $C_{19}$ and $C_{31}$ are nonhyperelliptic.
\end{lemma}

\begin{proof}
Let $P=(0:-1:1)$ on $F_p$. On $Z=1$ the function $X$ is a local
parameter at $P$, and $dX/(pY^{p-1})$ has order zero there.
A monomial $X^iY^jZ^l$ of degree $p-3$ gives a residue of
$h$-weight $j+1-k(i+1)$ modulo $p$.

For $p=19$ take $k=7$. Let $Q_1,Q_2,Q_3$ be the cyclic sums of
the monomials with exponent triples
\[
 (0,6,10),\qquad (1,13,2),\qquad (3,8,5),
\]
respectively. For $p=31$ take $k=5$ and use the two triples
\[
 (0,4,24),\qquad (7,8,13).
\]
All the monomials and their cyclic rotations have weight zero.
Their cyclic sums therefore give differentials on $C_p$.
In each sum the least exponent of $X$ occurs once. Its coefficient
at $Y=-1$ is nonzero. The orders at $P$ are thus $0,1,3$ for $p=19$
and $0,7$ for $p=31$. In each case the first order is zero, so the
map $F_p\to C_p$ is unramified at $P$. The orders downstairs are
the same.

On a hyperelliptic curve of genus $g$, the canonical vanishing
sequence is $0,1,\dots,g-1$ at an unramified point of
the hyperelliptic map and $0,2,\dots,2g-2$ at a ramification point.
Indeed, the canonical series is the pullback of
$|\mathcal O_{\mathbb P^1}(g-1)|$. For $g=3$ neither sequence contains $3$. For $g=5$ neither contains
$7$. This proves the lemma.
\end{proof}

\begin{remark}\label{rem:X15}
For $p=19$ the three forms are independent by their orders at $P$.
They form a basis since $g(C_{19})=3$. By Lemma~\ref{lem:small},
$C_{19}$ is nonhyperelliptic. Its canonical image in the coordinates
$[Q_2:-Q_1:Q_3]$ is therefore an irreducible plane quartic over $\Q$.
Let $P_{15}$ be the equation for $X_{15}$ in the introduction to
\cite{KLLRSS}. On $Z=1$, the polynomial $P_{15}(Q_2,-Q_1,Q_3)$
has zero remainder modulo $X^{19}+Y^{19}+1$. Thus $P_{15}$ vanishes
on the canonical image. Since it has degree four, it defines that
image. This gives $C_{19}\cong_{\Q}X_{15}$.
\end{remark}

\begin{proof}[Proof of Theorem~\ref{thm:B} and Corollary~\ref{cor:C}]
The curves descend to $\Q$ by Remark~\ref{rem:descent}.
Proposition~\ref{prop:toric-gonality} gives \eqref{eq:main-gonality}
and proves the automorphism assertion for $p>36$.
For $p=19,31$ apply Lemma~\ref{lem:small},
Proposition~\ref{prop:knownCM} and Lemmas~\ref{lem:roots-proper}
and~\ref{lem:aut-field}. The genus-one curve $C_7$ admits translations
over $\C$. The genus-two curve $C_{13}$ has a hyperelliptic involution.
This proves Theorem~\ref{thm:B}.

Finally, there are arbitrarily large primes $p\equiv1\pmod3$.
Since $g(C_p)=(p-1)/6$ and $\gon(C_p)\ge\sqrt p/3$, both the genus
and the gonality tend to infinity along these primes. Proposition~\ref{prop:knownCM} gives absolute simplicity
and complex multiplication.
\end{proof}

The curves $C_{31}$ and $D_{11}$ both have genus five and trivial
geometric automorphism group. Their Jacobians are absolutely simple
and are not isogenous over $\C$. Indeed, their endomorphism algebras
are $K_{31}$ and $\Q(\zeta_{11})$. The first has only the roots of
unity $\pm1$ and the second contains $\zeta_{11}$.

\section*{Acknowledgements}

The author acknowledges partial financial support from CNPq,
grant no.~302774/2025-4, and FAPESP, grant no.~2024/00923-6.

\section*{AI Declaration}

The author used ChatGPT (OpenAI, GPT-6 Astra Pro) to help extend the
proof to composite degrees and construct the genus-five isogeny.
The author also used it to locate references, revise the English prose,
and check calculations while reviewing the manuscript. The author takes
full responsibility for the final text and all mathematical results.

\end{document}